\documentclass[12pt]{amsart}
\usepackage{graphicx}
\usepackage{verbatim}
\usepackage{amsmath,amssymb}

\newtheorem{theorem}{Theorem}[section]
\newtheorem{lemma}[theorem]{Lemma}
\newtheorem{proposition}[theorem]{Proposition}
\newtheorem{corollary}[theorem]{Corollary}

\theoremstyle{definition}

\theoremstyle{remark}
\newtheorem{remark}[theorem]{Remark}
\newtheorem*{ack}{Acknowledgments}
\theoremstyle{definition}
 \newtheorem{definition}[theorem]{Definition}

\newcommand{\sign}{\operatorname{sign}}
\newcommand{\lk}{\operatorname{lk}}

\begin{document}

\title{Concordance Crosscap Number of A Knot}
\author{Gengyu Zhang}
\address{Department of Mathematics, Tokyo Institute of Technology, Tokyo 152-8551, Japan}
\email{zhang@math.titech.ac.jp} 
\date{\today}
\thanks{The author is supported by a Research Fellowship of the Japan Society for the Promotion of Science for Young Scientists.
This work is partly supported by the Grant-in-Aid for JSPS Fellows, the Ministry of Education, Culture, Sports, Science and Technology, Japan.
}

\begin{abstract}
We define the concordance crosscap number of a knot as 
the minimum crosscap number among all the knots concordant to the 
knot. The four-dimensional crosscap number is the minimum first 
Betti number of non-orientable surfaces smoothly embedded in 4-dimensional 
ball, bounding the knot. Clearly the 4-dimensional crosscap 
number is smaller than or equal to the concordance crosscap 
number. We construct two infinite sequences of knots to explain 
the gap between the two. In particular, the knot $7_4$ is one 
of the examples. 

\end{abstract}


\maketitle
\section{Introduction}
We work in the smooth category. Let $K$ be a knot embedded in the 3-sphere $S^3$. It can bound various orientable or non-orientable surfaces in $S^3$ or in the four-ball $B^4$ with $\partial B^4=S^3$. 

For orientable surfaces, one studies the genus.  The \textit{genus} $g(K)$ of a knot $K$ is the minimum
genus among all the orientable surfaces that $K$ bounds in $S^3$ and 
the \textit{slice genus} $g^*(K)$ is the minimum genus of orientable surfaces that $K$ bounds in $B^4$.
From the viewpoint of the concordance group, one can define the concordance genus $g_c(K)$ as the minimum genus in the concordance class of $K$. 
It is easy to see that $g^*(K)\leq g_c(K) \leq g(K)$. In response to the
question asked by Gordon \cite{Gordon:problem}, Nakanishi explained the gap between $g^*(K)$ and $g_c(K)$ in \cite{Nakanishi:unknotting}. The relation among these knot invariants are also investigated by Livingston in \cite{Livingston:concordance_genus}.
  
For non-orientable surfaces, instead of the genus, the first Betti number is taken as an invariant. 
The \textit{crosscap number} $\gamma(K)$, also known as \textit{non-orientable genus},  
of $K$ is defined to be the minimum first Betti number
of non-orientable surfaces in $S^3$ bounding $K$. See \cite{Clark:crosscap, Murakami/Yasuhara:crosscap, Teragaito:crosscap,
Hirasawa/Teragaito:crosscap, Ichihara/Mizushima:crosscap} for studies on this invariant. 
The \textit{4-dimensional crosscap number} $\gamma^*(K)$ is the minimum first 
Betti number of non-orientable surfaces in the 4-ball $B^4$ bounding $K$. Some results on 4-dimensional crosscap number 
by Murakami and Yasuhara can be found in \cite{Murakami/Yasuhara:4genus}. 
In this paper we define the \textit{concordance crosscap
number} $\gamma_c(K)$ as the minimum crosscap number among all the knots concordant to $K$, and study the gap between
$\gamma^*(K)$ and $\gamma_c(K)$. 

We give a necessary condition for pretzel knots of type $P(4-p, p, 2n-p)$ and $P(-1-p, p,2n-p)$ to have concordant crosscap number 1, 
and construct a series of knots such that there exist infinitely many knots 
with $\gamma^*(K)< \gamma_c(K)$. In particular, the knot $7_4$ is one of the examples. 

\section{Preliminaries}

Two knots $K_0$ and $K_1$ in the 3-sphere $S^3$ are called $concordant$ if they cobound an annulus in $S^3\times [0,1]$ 
with $K_0$ and $K_1$ in $S^3\times \{0\}$ and $S^3\times \{1\}$ respectively. Concordance is an equivalence relation and the set
of the equivalence classes forms an abelian group, the so called \textit{concordance group}. The \textit{concordance crosscap
number} is defined from this point of view.

Given a Seifert surface $F$ bounding $K$, one can associate a Seifert matrix $V$ to it. The \textit{signature} $\sigma(K)$ of $K$ is defined to be
the signature of the symmetric matrix $V+V^T$, i.e. $\sigma(K)=\sign(V+V^T)$. It is well known that the knot signature is
a concordance invariant \cite{Murasugi:numerical_invariant}. 

Given a surface bounding a knot, not necessarily orientable, one can associate a Goeritz matrix with it, see for example \cite{Lickorish:knot}.
Gordon and Litherland in \cite{Gordon/Litherland:signature} gave a simple algorithm to calculate the
knot signature, expressed as the signature of the Goeritz matrix corresponding to the knot plus a correction term. 

\begin{lemma}[\cite {Gordon/Litherland:signature}\label{thm:signature}]
Let $F$ be any surface bounding a knot $K$. Then the signature $\sigma (K)$ can be calculated out of the Goeritz matrix $G_F$ 
and the normal Euler number $e(F)$, namely we have
\begin{equation*}
 \sigma (K)=\sign(G_F) +\dfrac{1}{2}e(F),
\end{equation*} 
where $e(F) :=-\lk(K, K')$. Here $K':=\partial N(K) \cap F$ and $N(K)$ denotes the regular neighborhood of the knot $K$. 
\end{lemma}

Note that the normal Euler number vanishes for an orientable surface $F$. In this case, the result conforms
to the original definition. 

In \cite{Clark:crosscap} Clark introduced the crosscap number, a knot invariant from
a non-orientable viewpoint. The following result is well known.   
\begin{lemma}[\cite{Clark:crosscap}] \label{thm:crosscap1}
  A knot $K$ has crosscap number $1$ if and only if it is a $({2,n})$-cable knot. 
\end{lemma}
\begin{definition}
The \textit{concordance crosscap number} of a knot $K$, denoted by $\gamma_c(K)$, is the minimum crosscap number among all the knots in the same concordance
class as $K$.  
\end{definition}
The \textit{4-dimensional crosscap number}, or \textit{non-orientable 4-genus} \cite{Murakami/Yasuhara:4genus} $\gamma^*(K)$ is 
the minimum first Betti number of non-orientable surfaces in $B^4$ bounding the knot $K$. The concordance
crosscap number can be regarded as a bridge to relate the 4-dimensional invariant and the 3-dimensional one. 
\begin{remark}
By convention, the crosscap number of the unknot is defined to be 0 for completeness. 
So if the knot $K$ is a slice knot, we define $\gamma^*(K)=\gamma_c(K)=0$. 
\end{remark}

Clearly  $\gamma^*(K)\leq \gamma_c(K)\leq \gamma(K)$, and the inequalities are best possible. 
 
\begin{proposition} \label{thm:cable}
Let $K$ be a $(2, p)$-cable knot with $p\neq \pm 1$. Then $\gamma^*(K)= \gamma_c(K)= \gamma(K)=1$. 
\end{proposition}
\begin{proof}
It is known that $\gamma(K)=1$ by Lemma~\ref{thm:crosscap1}, and we will prove that $\gamma^\ast(K)\neq 0$.
Suppose that $K$ is the $(2,p)$-cable of a knot, then 
it bounds a twisted M\"obius band, whose corresponding Goeritz matrix is the 1 by 1 matrix $\begin{pmatrix} p \end{pmatrix}$,
 and whose normal Euler number is $-2p$. By Lemma~\ref{thm:signature} we have $\sigma(K)=\sign(p)-p$. 
Then $\sigma(K)\neq 0$ since $p\neq \pm 1$, and $K$ is not a slice knot. Then $\gamma^\ast(K)\geq 1$ and the proof is finished.  
\end{proof}
Then we ask whether some knot exists such that all knots in its concordance class have crosscap numbers bigger than the 
4-dimensional crosscap number of $K$, i.e. whether there exists a knot $K$ with $\gamma^*(K)< \gamma_c(K)$.






\section{Examples constructed from Seifert surfaces}
\begin{figure}
  \includegraphics [scale=.38] {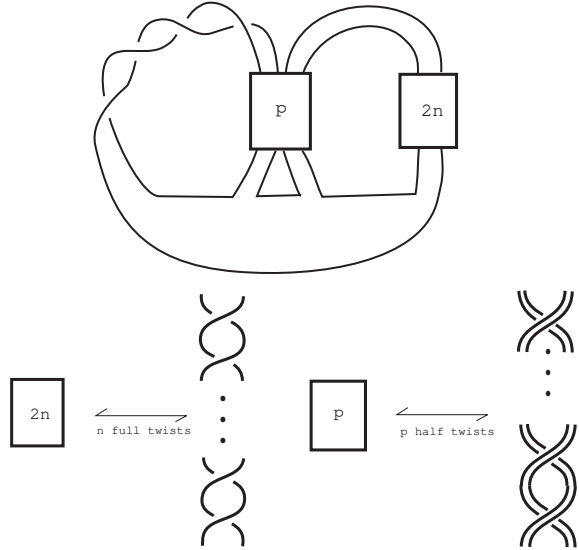}
  \caption{$K(4, 2n, p)$}
  \label{fig:orientable}
\end{figure}

Consider the knot illustrated in Figure~\ref{fig:orientable}. We denote it by $K(4, 2n, p)$, where $p$ is an odd number and $n\in \mathbf{Z}$ 
denotes the number of the full twists on the right band.  It is easy to see that the knot diagram actually
represents the pretzel knot $P(4-p, p, 2n-p)$. 

\begin{figure}
  \includegraphics [scale=.4] {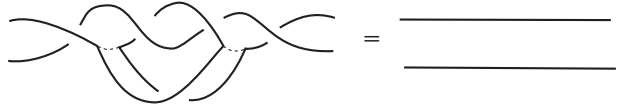}
  \caption{Untwisted by adding a twisted band}
  \label{fig:unknot}
\end{figure}

By adding a twisted band as illustrated in Figure~\ref{fig:unknot}, the knot $K(4, 2n, p)$ is changed into a ribbon knot. 
To see this, consider the indicated fission band in Figure~\ref{fig:slice}, which gives an unlink. Reversing the fission 
process, one will see a fusion band intersecting two disks bounding the unlink in arcs.  
Since the ribbon knot is a slice knot \cite{Rolfsen:1990},   
the disk in $B^4$ that it bounds and the twisted band together give a first Betti number 1 non-orientable surface in $B^4$
bounding the knot $K(4, 2n, p)$, and we have $\gamma^*(K(4, 2n, p))\leq 1$ immediately.

\begin{theorem}\label {thm:orientable}
If $\gamma_c(K(4, 2n, p))=1$, there exists some odd number $l$ such that the following equalities hold. 
\begin{equation*}
  \begin{cases}
   p^2-8n=l^2  &\text{when} \mspace{5mu}8n-p^2<0,\\
   3(8n-p^2)=l^2 &\text {when} \mspace{5mu}8n-p^2>0. 
  \end{cases}
\end{equation*}
\end{theorem}

\begin{proof} 

The Goeritz matrix $G_{F_{p,n}}$ of the knot $K(4, 2n, p)$ with respect to the Seifert surface $F_{p,n}$ under a generator system \{$a,b$\}, 
indicated with oriented dashed lines in Figure~\ref{fig:generator},
of the homology group $H_1(F_{p,n};\mathbf{Z})$ is of the form $\begin{pmatrix}4&-p\\-p&2n \end{pmatrix}$. Since the surface $F_{p,n}$ 
is orientable, then the normal Euler number $e(F_{p,n})$ is $0$ and the signature of $K(4, 2n, p)$ is just the signature of the Goeritz matrix $G_{F_{p,n}}$ according to Lemma~\ref{thm:signature}.

\begin{figure}
\includegraphics [scale=.38] {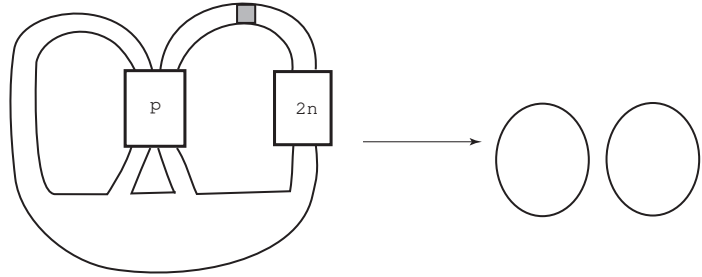}
\caption{Fission band giving an unlink}
\label{fig:slice}
\end{figure}

\begin{figure}
  \includegraphics [scale=.38] {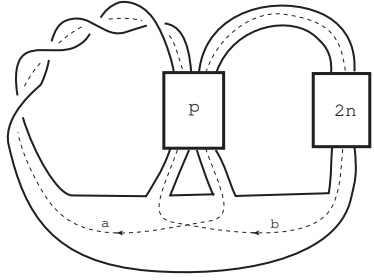}
  \caption{Generator system \{$a, b$\}}
  \label{fig:generator}
\end{figure}

Note that the value of $\sign G_{F_{p,n}}$
is determined by the sign of the $2\times 2$ matrix's determinant, i.e., $\sign G_{F_{p,n}}=0$ when $\det G_{F_{p,n}}<0$, 
and $\sign G_{F_{p,n}}=2$ when $\det G_{F_{p,n}}>0$.

When $\gamma_c(K(4, 2n, p))=1$,  there exists a knot $K'$ concordant to $K(4, 2n, p)$ with $\gamma(K')=1$. 
So $K'$ is a 2-cable of some knot by Lemma~\ref{thm:crosscap1}, and $\sigma(K')=\sign(q)-q$ by the proof of Proposition~\ref{thm:cable}. 
Due to the fact that the signature is a concordance invariant, we also have $\sigma(K(4, 2n, p))=\sigma(K')$.
Then by the above discussion $\sigma(K')=\sign(q)-q \in \{2, 0\}$ and a calculation we obtain the following values:
\begin{equation*}
  \begin{cases}
   |q|=1 &\text{when} \mspace{5mu}\det G_{F_{p,n}}<0,\\
   q=3 &\text{when} \mspace{5mu}\det G_{F_{p,n}}>0.
  \end{cases}
\end{equation*}

On the other hand, by concordance it is well known that the knot $K(4, 2n, p)\# (-\overline {K'})$ is a slice knot where $-\overline {K'}$ denotes the
mirror image of $K'$ with a reverse orientation. 
Then we have the Alexander polynomial $\Delta(K(4, 2n, p)\# (-\overline{K'}))=\Delta(K, 2n, p)\Delta(K')=f(t)f(t^{-1})$
for some polynomial $f(t)$ with $f(1)=1$ (see \cite{Fox/Milnor:cobordism}). It follows that $\det(K')\det(K(4, 2n, p))=l^2$ for some odd number $l$,
where $\det(K)$, the \textit{determinant} of the knot $K$, is equal to $|\Delta_K(-1)|$. 
It is also known to be the order of the first homology group of double branched cover of $S^3$ branched 
over $K$, and so we have $\det(K(4, 2n, p))=|\det G_{F_{p,n}}|=|8n-p^2|$. Thus $|8n-p^2|\cdot |q|=l^2$, and precisely speaking 
\begin{equation*}
  \begin{cases}
   p^2-8n=l^2  &\text{when} \mspace{5mu}8n-p^2<0,\\
   3(8n-p^2)=l^2 &\text {when} \mspace{5mu}8n-p^2>0. 
  \end{cases}
\end{equation*}
\end{proof}

\begin{corollary}\label{thm:ori}
 For any odd number $p$, there exist infinitely  many $n$ 
 for which the concordance crosscap number of the knot $K(4, 2n, p)$ is greater than $1$ and 
4-dimensional crosscap number equal to $1$. 
\end{corollary}
\begin{proof}
In the case of $8n-p^2>0$, we have $\sigma(K(4, 2n, p))=2$, and so the knot $K(4, 2n, p)$ is not slice. Then
we have $\gamma^*(K(4,2n,p))=1$, and $\gamma_c(K(4,2n,p))\geq 1$. 
By Theorem~\ref{thm:orientable} for a knot with $\gamma_c(K(4,p,2n))=1$ and $8n-p^2>0$ the equality $3(8n-p^2)=l^2$ holds. 
Then the number $3$ either becomes a common divisor of $8n$ and $p$, or divides neither $8n$ nor $p$.  While there exist infinitely many $n$ 
with $3$ dividing only one of $8n$ and $p$, in which cases $3(8n-p^2)$ cannot make a square for those $n$, 
therefore a contradiction arises. The result follows. 

\end{proof}

Therefore Corollary~\ref{thm:ori} gives examples
demonstrating the gap between the 4-dimensional crosscap number and the concordance crosscap number. 
One can also discuss the knot of type $K(-4, 2n, p)$ 
with reverse twists on the left band in the same way. 
Note that the crosscap number of the pretzel knot $P(i_1,i_2, \cdots, i_n)$ with all $p_i$ odd is shown to be $n$ in \cite{Ichihara/Mizushima:crosscap}.
Then we know $\gamma(K(4,2n,p))=\gamma(K(-4,2n,p))=3$, but we cannot say whether the concordance crosscap number of knot $K(4,2n,p)$(or $K(-4,2n,p)$) 
is 3 or not.

\begin{remark}
The knot $K(-4,-4,-1)=P(-3,-1,-3)$ is knot $7_4$ in Alexander-Briggs' notation \cite{Alexander/Briggs:knotted, Rolfsen:1990}.
\end{remark}

\begin{figure}
  \includegraphics [scale=.35] {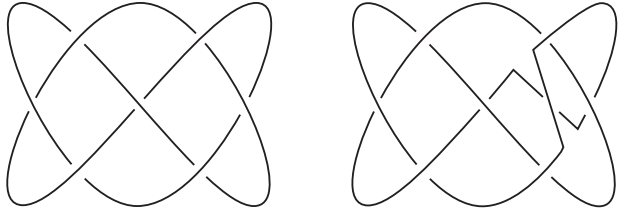}
  \caption{Knot $7_4$ and unknotting of $7_4$ by adding a twisted band}
  \label{fig:knot74}
\end{figure}

\section{Examples constructed from non-orientable surfaces}
\begin{figure}
  \includegraphics [scale=.28] {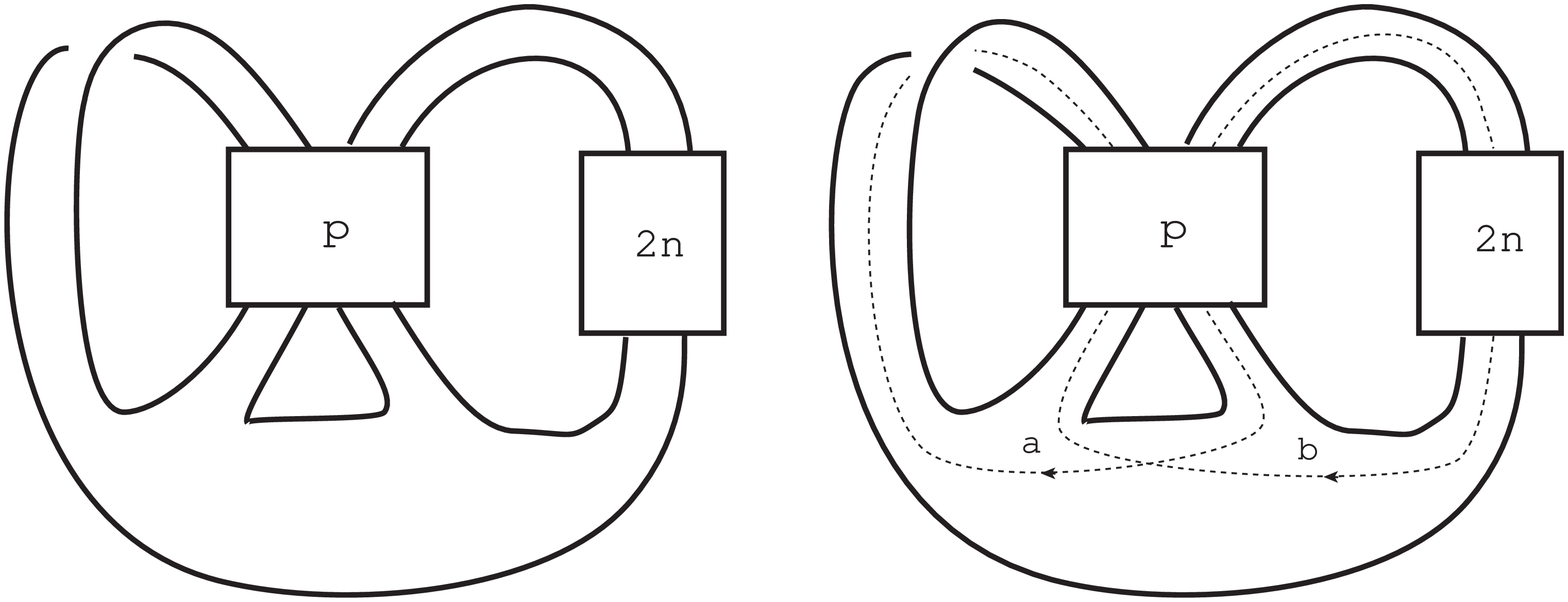}
  \caption{$K(-1, 2n, p)$}
  \label{fig:nonorientable}
\end{figure}

\begin{figure}
  \includegraphics [scale=.46] {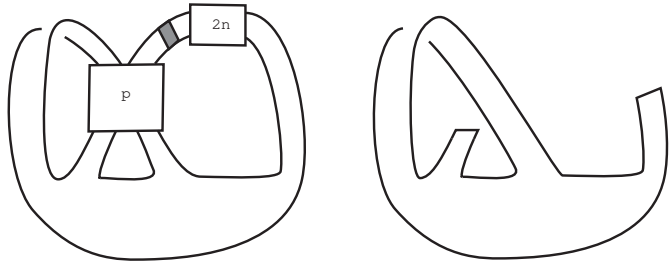}
  \caption{Unknotting by adding a twisted band}
  \label{fig:unknot2}
\end{figure}

Consider the knot $K(-1, 2n, p)$ illustrated in Figure~\ref{fig:nonorientable}, which is a pretzel knot of type $P(-1-p,p,2n-p)$. 
 Here the $n$ and $p$ are the same as in Section~3, but in order to make our statement clearer and more succinct
we only discuss the case with $n>0$. The case $K(1, 2n, p)$ with $n<0$
can be realized by taking a mirror image, and discussed in the same way. It is obvious that $\gamma^*(K(-1,2n,p))\leq 1$ since it can be unknotted by adding
a half twisted band \cite{Lickorish:unknotting}. See Figure~\ref{fig:unknot2}.  

Note that the knot $K(-1, 2n, \pm 1)$ is actually 
the torus knot $T(2, 2n+1)$ with $\gamma^*=\gamma_c=\gamma=1$ by Proposition~\ref{thm:cable}.

\begin{theorem}\label{thm:non}
If $\gamma_c(K(-1, 2n, p))=1$, there exists an odd number $l$ such that $(2n+p^2)(2n+1)=l^2$ when $n>0$.
\end{theorem}

\begin{proof}
The knot diagram bounds a non-orientable surface $F'_{p,n}$ with normal Euler number $-4n$. The Goeritz matrix $G_{F'_{p,n}}$ 
corresponding to the non-orientable surface with the indicated basis \{$a, b$\} (see Figure~\ref{fig:nonorientable}) of the homology group $H_1(F'_{p,n})$
is of the form $\begin{pmatrix} -1&-p\\-p&2n\end{pmatrix}$. Thus with $n>0$ we have $\det G_{F'_{p,n}}=-2n-p^2<0$, and hence
$\det K(-1, 2n, p)=2n+p^2$. By Lemma~\ref{thm:signature} $\sigma(K(-1, 2n, p))=-2n$ since 
$\sign G_{F'_{p,n}}=0$ due to the negative sign of the matrix determinant. 
 
If $\gamma_c(K(-1, 2n, p))=1$, through a similar argument as the knot of type $K(4, 2n, p)$, 
it is concordant to the ($2, 2m+1$)-cable of a knot $K'$, with signature $\sigma(K')=\sign(2m+1)-(2m+1)$
and determinant $\det(K')=2m+1$. By concordance
we have $-2n=\sign(2m+1)-(2m+1)$ which gives us the equality $n=m$ by a calculation. 
Then the determinant of $K'$ becomes $2n+1$, and therefore
there exists some odd number $l$ such that $\det K' \det K(-1,2n,p)=(2n+1)(2n+p^2)=l^2$ holds. 

\end{proof}

\begin{remark}
Through the proof of Theorem~\ref{thm:non} we know that the signature of the knot $K(-1,2n,p)$ is $-2n$, 
and so different value of $n$ gives us different concordance class of the knot type $K(-1, 2n, p)$.
Therefore we can cover as many concordance classes as necessary. 
\end{remark}

\begin{corollary}\label{thm:nonori}
For any odd number $p\neq \pm 1$ there exist infinitely many $n$ 
for which the concordance crosscap number of the knot $K(-1, 2n, p)$ is greater than $1$ and 4-dimensioanl crosscap
number equal to $1$.
\end{corollary}
\begin{proof}
By Theorem~\ref{thm:non} with $n>0$ we have $(2n+p^2)(2n+1)=l^2$ if $\gamma_c(K(-1, 2n, p))=1$. 
In this case $\sigma(K(-1,2n,p))<0$, thus the knot $K(-1,2n,p)$ is not slice and $\gamma^*(K(-1,2n,p))=1$. 
It is well known that infinitely many primes are in the set $\{2m+1|2m+1>p^2 \}=\{2n+p^2|n\in \mathbf{Z}^+\}$. 
 It is clear that $2n+1$ cannot be divided by $2n+p^2$, which confirms that those $(2n+p^2)(2n+1)$ with prime $2n+p^2$ cannot make a square. 
Therefore for the knots $K(-1,2n,p)$ with those $n$ the concordance crosscap numbers are greater than 1. Then we have the proof. 
\end{proof}

\begin{remark}
In \cite{Ichihara/Mizushima:crosscap} the pretzel knot $P(i_1, i_2, \cdots, i_n)$ with one $i_k$ even and other odd 
is shown to have crosscap number $n-1$. So for the knots in Corollary~\ref{thm:nonori}, we have $\gamma_c(K(-1, 2n, p))=2$ since $\gamma(P(-1-p, p, 2n-p))=2$ for any
($p, n$). 
\end{remark}

In conclusion there are also numerous knots of type $K(-1, 2n, p)$ with concordance crosscap number larger than their 4-dimensional crosscap
number.

\begin{ack}
The author would like to express her thanks to  Professor Hitoshi Murakami and Professor Akira Yasuhara for helpful 
discussions and encouragement while pursuing this topic. She also would like to thank the participants in the Workshop on Crosscap
Number, held in Shikotsuko Hokkaido, Japan between August 1st and 4th 2006, for pointing out an error in the previous version. 
\end{ack}

\bibliography{mrabbrev,gengyu}

\providecommand{\bysame}{\leavevmode\hbox to3em{\hrulefill}\thinspace}
\providecommand{\MR}{\relax\ifhmode\unskip\space\fi MR }
\providecommand{\MRhref}[2]{%
  \href{http://www.ams.org/mathscinet-getitem?mr=#1}{#2}
}
\providecommand{\href}[2]{#2}
\begin{thebibliography}{10}

\bibitem{Alexander/Briggs:knotted}
J.~W. Alexander and G.~B. Briggs, \emph{On types of knotted curves}, Ann. of
  Math. (2) \textbf{28} (1926/27), no.~1-4, 562--586. \MR{MR1502807}

\bibitem{Clark:crosscap}
B.~E. Clark, \emph{Crosscaps and knots}, Internat. J. Math. Math. Sci.
  \textbf{1} (1978), no.~1, 113--123. \MR{0478131 (57 \#17620)}

\bibitem{Fox/Milnor:cobordism}
R.~H. Fox and J.~W. Milnor, \emph{Singularities of {$2$}-spheres in {$4$}-space
  and cobordism of knots}, Osaka J. Math. \textbf{3} (1966), 257--267.
  \MR{0211392 (35 \#2273)}

\bibitem{Gordon:problem}
C.~McA. Gordon, \emph{Problems}, Knot theory (Proc. Sem., Plans-sur-Bex, 1977),
  Lecture Notes in Math., vol. 685, Springer, Berlin, 1978, pp.~309--311.

\bibitem{Gordon/Litherland:signature}
C.~McA. Gordon and R.~A. Litherland, \emph{On the signature of a link}, Invent.
  Math. \textbf{47} (1978), no.~1, 53--69. \MR{0500905 (58 \#18407)}

\bibitem{Hirasawa/Teragaito:crosscap}
M.~Hirasawa and M.~Teragaito, \emph{{Crosscap numbers of 2-bridge knots}},
  Topology \textbf{45} (2006), no.~3, 513--530.

\bibitem{Ichihara/Mizushima:crosscap}
K.~Ichihara and S.~Mizushima, \emph{Crosscap numbers of pretzel knots},
  preprint, 2006.

\bibitem{Lickorish:unknotting}
W.~B.~R. Lickorish, \emph{Unknotting by adding a twisted band}, Bull. London
  Math. Soc. \textbf{18} (1986), no.~6, 613--615. \MR{859958 (88e:57008)}

\bibitem{Lickorish:knot}
\bysame, \emph{An {I}ntroduction to {K}not {T}heory}, Graduate Texts in
  Mathematics, vol. 175, Springer-Verlag, New York, 1997. \MR{1472978
  (98f:57015)}

\bibitem{Livingston:concordance_genus}
C.~Livingston, \emph{The concordance genus of knots}, Algebr. Geom. Topol.
  \textbf{4} (2004), 1--22 (electronic). \MR{2031909 (2005e:57023)}

\bibitem{Murakami/Yasuhara:crosscap}
H.~Murakami and A.~Yasuhara, \emph{Crosscap number of a knot}, Pacific J. Math.
  \textbf{171} (1995), no.~1, 261--273. \MR{1362987 (96k:57006)}

\bibitem{Murakami/Yasuhara:4genus}
\bysame, \emph{Four-genus and four-dimensional clasp number of a knot}, Proc.
  Amer. Math. Soc. \textbf{128} (2000), no.~12, 3693--3699. \MR{1690998
  (2001b:57020)}

\bibitem{Murasugi:numerical_invariant}
K.~Murasugi, \emph{On a certain numerical invariant of link types}, Trans.
  Amer. Math. Soc. \textbf{117} (1965), 387--422. \MR{0171275 (30 \#1506)}

\bibitem{Nakanishi:unknotting}
Y.~Nakanishi, \emph{A note on unknotting number}, Math. Sem. Notes Kobe Univ.
  \textbf{9} (1981), no.~1, 99--108. \MR{634000 (83d:57005)}

\bibitem{Rolfsen:1990}
D.~Rolfsen, \emph{Knots and {L}inks}, Mathematics Lecture Series, vol.~7,
  Publish or Perish Inc., Houston, TX, 1990. \MR{95c:57018}

\bibitem{Teragaito:crosscap}
M.~Teragaito, \emph{Crosscap numbers of torus knots}, Topology Appl.
  \textbf{138} (2004), no.~1-3, 219--238. \MR{2035482 (2005a:57007)}

\end{thebibliography}
\bibliographystyle{amsplain}

\end{document}